\documentclass[10pt, letterpaper]{amsart}
\usepackage[latin1]{inputenc}
\usepackage{amsfonts}
\usepackage{amsmath}
\usepackage{amssymb}
\usepackage{amsthm}
\usepackage[all]{xy}
\usepackage{graphicx}
\usepackage[margin=1.25in]{geometry}

\theoremstyle{plain}
	\newtheorem*{theorem*}{Theorem}
	\newtheorem{theorem}{Theorem}[subsection]
	\newtheorem*{lemma*}{Lemma}
	\newtheorem{lemma}{Lemma}[subsection]
	\newtheorem*{claim*}{Claim}
	
	\newtheorem*{prop*}{Proposition}
	\newtheorem{prop}{Proposition}[subsection]
	
	\newtheorem*{cor*}{Corollary}

\theoremstyle{definition}
	\newtheorem{note}{Remark}[subsection]
	\newtheorem{definition}{Definition}[subsection]
	
\newcommand{\Z}{\mathbb{Z}}
\newcommand{\R}{\mathbb{R}}

\newcommand{\C}{\mathbb{C}}
\newcommand{\A}{\mathbb{A}}

\newcommand{\isoarrow}{\ensuremath{\overset{\sim}{\longrightarrow}}}
\newcommand{\noqed}{\renewcommand{\qedsymbol}{}}

\begin{document}

\title[An exact formula: lattice points and the automorphic spectrum]{An exact formula relating lattice points in symmetric spaces to the automorphic spectrum}
\author{Amy T. DeCelles}
\address{University of Minnesota \\ School of Mathematics \\127 Vincent Hall \\206 Church St. SE \\Minneapolis, Minnesota 55455}
\email{decel004@umn.edu}
\urladdr{www.math.umn.edu/~decel004}
\thanks{This paper presents results from the author's PhD thesis, completed under the supervision of Professor Paul Garrett, whom the author thanks warmly for many helpful conversations.  The author was partially supported by the Doctoral Dissertation Fellowship from the Graduate School of the University of Minnesota and by NSF grant DMS-0652488.}


\begin{abstract}We extract an exact formula relating the number of lattice points in an expanding region of a complex semi-simple symmetric space and the automorphic spectrum from a spectral identity, which is obtained by producing two expressions for the automorphic fundamental solution of the invariant differential operator $(\Delta - \lambda_z)^{\nu}$.  On one hand, we form a Poincar\'{e} series from the solution to the corresponding differential equation on the free space $G/K$, which is obtained using the harmonic analysis of bi-$K$-invariant functions.  On the other hand, a suitable global automorphic Sobolev theory, developed in this paper, enables us to use the harmonic analysis of automorphic forms to produce a solution in terms of the automorphic spectrum.  \end{abstract}

\maketitle

\section{Introduction}

The simplest lattice-point counting problem is the \emph{Gauss circle problem}, counting lattice points within a circle in the Euclidean plane.  Elementary packing arguments yield
$$N(T) \;\;=  \;\; \# \{ \xi \in \Z^2 : | \xi | \leq T \} \;\; = \;\; 2\pi \cdot T^2 \; + \; \operatorname{O}(T)$$
Determining the optimal error term, conjectured to be $\operatorname{O}(T^{1/2 + \varepsilon})$, is a topic of active research \cite{iwaniec-mozzochi1988,huxley2003}.  
In the hyperbolic plane, where the circumference of a circle is proportionate to its area, packing arguments fail to produce an asymptotic with error term.  Subtler methods have produced asymptotics for lattice-point counting in hyperbolic spaces; see \cite{patterson1975, miatello-wallach1992,bruggeman-miatello-wallach1999}.   In affine symmetric spaces, ergodic methods have produced asymptotics for lattice point counting, e.g. \cite{bartels1982,duke-rudnick-sarnak1993,eskin-mcmullen1993,maucourant2007}.  Gorodnik and Nevo \cite{gorodnik-nevo2010} give a good exposition of ergodic methods and lattice-point counting.

In contrast, spectral methods produce an \emph{exact formula} relating the number of lattice points in an expanding region and the automorphic spectrum.  Instead of expressing something mysterious in terms of something familiar, this formula, like the explicit formula of Riemann-von Mangoldt, which relates the prime numbers to the zeros of zeta, gives a relationship between two mysterious things; its appeal lies not in its utility for evaluation of one side, but rather in that it reveals a connection between seemingly disparate things.


Instead of using trace formula methods, we obtain a spectral identity from a \emph{differential equation} on the arithmetic quotient $X = \Gamma \backslash G/K$.  This offers the possibility of operating under somewhat different hypotheses than one would usually take when using trace formula or relative trace formula methods: while trace formulas work best with very smooth data, the applications to asymptotics of L-functions that we have in mind (see \cite{decellesPhD2011}) use a Poincar\'{e} series whose data is not smooth nor compactly supported.  On the other hand, our approach necessitates a more careful treatment of analytic issues.  There have been relatively few papers discussing the analytic issues involved in the relative trace formula; see \cite{jacquet1995,finis-lapid-mueller2009}.

We obtain a spectral identity by producing two different expressions for the solution of the following differential equation on the arithmetic quotient $X = \Gamma \backslash G/K$:
$$(\Delta - \lambda_z)^{\nu} \; v_z \;\; = \;\; \delta_{x_o}$$
where the Laplacian $\Delta$ is the image of the Casimir operator for $\mathfrak{g}$, $\lambda_z$ a complex parameter, $\nu$ an integer, and $\delta_{x_o} = \delta_{\Gamma \cdot 1\cdot K}$ the Dirac delta distribution at the basepoint in $\Gamma \backslash G/K$.  On one hand, we form a Poincar\'{e} series from the solution to the corresponding differential equation on the free space $G/K$, obtained using the analysis of bi-$K$-invariant functions \cite{decelles2011}.  On the other hand, a global automorphic Sobolev theory produces a solution in terms of the automorphic spectrum.  For a sketch of this discussion for $SL_2(\C)/SU(2)$, see \cite{garrett-newark2010}.


\section{Global automorphic Sobolev theory}

\subsection{Parametrization of spectrum, spectral transform and inversion} \label{sp_param}

The spectral theory of automorphic forms decomposes square integrable automorphic forms in terms of eigenfunctions.  For succinctness, we restrict our attention to automorphic forms that are  \emph{spherical} at infinity and have \emph{trivial central character}.   Though the automorphic spectrum consists of disparate pieces (cusp forms, Eisenstein series, residues of Eisenstein series) it will be useful to have a uniform notation.  We posit a parameter space $\Xi$ with spectral (Plancherel) measure $d\xi$ and let $\Phi_{\xi}$ denote the elements of the spectrum.  The general spectral theory (see, for example, \cite{langlands76,moeglin-waldspurger95}) implies the following. 

For test functions $f$ on $G_k \backslash G_{\A}$, the spectral transform $\mathcal{F}: C_c^{\infty}(G_k \backslash G_{\A}) \to C^0(\Xi)$ by $f  \to \langle f, \Phi_{\xi} \rangle$, where $\langle \, , \, \rangle$ is the usual $L^2$ inner product, extends to an isometry
$$\mathcal{F}: L^2(Z_{\A} G_k \backslash G_{\A}/K_{\infty}) \;\; \isoarrow \;\; L^2(\Xi)$$
On compactly supported continuous functions, the inverse map is given by an integral formula, and there is an inversion formula
$$f \;\; \overset{L^2}{=} \;\; \int_{\Xi} \mathcal{F}f(\xi) \, \Phi_{\xi} \, d\xi$$
converging at least for test functions on $Z_{\A}G_k \backslash G_{\A}$.

Since the critical issues for applications arise at the \emph{archimedean} place, we will consider global \emph{archimedean} spherical automorphic Sobolev spaces.  For a discussion of global automorphic Sobolev spaces for $SL_2(\C)$, see \cite{garrett2011-afcsob}.

\subsection{Characterizations of Sobolev spaces}

Let $G$ be a semi-simple or reductive Lie group with discrete subgroup $\Gamma$ and maximal compact subgroup $K$, and let $X$ be the arithmetic quotient $\Gamma \backslash G/K$.  We define positive index global archimedean spherical automorphic Sobolev spaces as right $K$-invariant subspaces of completions of $C_c^{\infty}(\Gamma \backslash G)$ with respect to a topology induced by seminorms associated to derivatives from the universal enveloping algebra, as follows.  Let $\mathcal{U}\mathfrak{g}^{\leq \ell}$ denote the finite dimensional subspace of the universal enveloping algebra $\mathcal{U}\mathfrak{g}$ consisting of elements of degree less than or equal to $\ell$.  Each $\alpha \in \mathcal{U}\mathfrak{g}$ gives a seminorm $\nu_{\alpha}(f) \; = \; \lVert \alpha f \rVert_{L^2(\Gamma \backslash G)}^2$  on $C_c^{\infty}(\Gamma \backslash G)$.

\begin{definition}Consider the space of smooth functions bounded with respect to these seminorms:
$$\{ f \in C^{\infty}(\Gamma \backslash G) : \nu_{\alpha} f \, < \, \infty\;\text{ for all } \alpha \in \mathcal{U}\mathfrak{g}^{\leq \ell}\}$$
Let $H^{\ell}(\Gamma \backslash G)$ be the completion of this space with respect to the topology induced by the family $\{\nu_{\alpha} : \alpha \in \mathcal{U}\mathfrak{g}^{\leq \ell} \}$.  The \emph{global archimedean spherical automorphic Sobolev space} $H^{\ell}(X) = H^{\ell}(\Gamma \backslash G)^K$ is the subspace of right-$K$-invariant functions in $H^{\ell}(\Gamma \backslash G)$.\end{definition}

\begin{prop} The space of test functions $C_c^{\infty}(X)$ is dense in $H^{\ell}(X)$.

\begin{proof}
We approximate a smooth function $f \in H^{\ell}(X)$ by pointwise products with smooth cut-off functions, as follows.  Let $\{\eta_n\}  \subset C_c^{\infty}(\Gamma \backslash G)$ be a family of left $N$-invariant, right $K$-invariant smooth cut-off functions with $\cup_n \mathrm{spt}(\eta_n) = \Gamma \backslash G$ and  $\sup_{g \in \Gamma \backslash G} | \alpha \eta_n (g) | \, \ll \, 1$ for all $\alpha \in \mathcal{U}\mathfrak{g}^{\leq \ell}$, where the implied constant does not depend on the support of $\eta$, but may depend on $\ell$.  By definition, $ \nu_{\gamma} \big(\eta_n \cdot f - f\big) \, = \, \lVert \gamma \big( \eta_n \cdot f - f\big) \rVert_{L^2(\Gamma \backslash G)}$.  Leibnitz' rule implies that $\gamma\big( \eta_n \cdot f - f\big)$ is a finite linear combination of terms of the form $\alpha(\eta_n - 1) \cdot \beta f$ where $\alpha$, $\beta \in \mathcal{U}\mathfrak{g}^{\leq \ell}$.  When $\alpha = 0$, 
$$\lVert \alpha(\eta_n -1)\cdot \beta f \rVert_{L^2(\Gamma \backslash G)}\;\; = \;\; \lVert (\eta_n -1) \cdot \beta f \rVert_{L^2(\Gamma \backslash G)} \;\; \leq \;\;  \int_{\Gamma \backslash G - \mathrm{spt} \eta_n} \; |(\beta \, f) (g)|^2 \, dg$$
Otherwise, $\alpha (\eta_n -1) = \alpha \eta_n$, and
\begin{eqnarray*}
\lefteqn{\lVert \alpha(\eta_n -1)\cdot \beta f \rVert_{L^2(\Gamma \backslash G)} \;\; = \;\; \lVert \alpha \eta_n \cdot \beta f \rVert_{L^2(\Gamma \backslash G)} }\\
 & \ll & \sup_{g \in G} \; |\alpha \, \eta_n(g) | \; \cdot  \; \int_{\Gamma \backslash G - \mathrm{spt} \eta_n} \; |(\beta \, f) (g)|^2 \, dg \;\; \ll\;\;  \int_{ \Gamma \backslash G - \mathrm{spt} \eta_n} \; |(\beta \, f) (g)|^2 \, dg
\end{eqnarray*}
Let $B$ be any bounded set containing all of the (finitely many) $\beta$ that appear as a result of applying Leibniz' rule.  Then
$$ \nu_{\gamma} \big(\eta_n \cdot f - f\big) \;\; \ll \;\; \sup_{\beta \in B} \;\;  \int_{\Gamma \backslash G - \mathrm{spt} \eta_n} \; |(\beta \, f) (g)|^2 \, dg$$
Since $B$ is bounded and $f \in H^{\ell}(X)$, the right hand side approaches zero as $ n \to \infty$.
\end{proof}
\end{prop}

\begin{prop} \label{afc_casimir_norm} Let $\Omega$ be the Casimir operator in the center of $\mathcal{U}\mathfrak{g}$.  The norm $\lVert \,\cdot \, \rVert_{2\ell}$ on $C_c^{\infty}(\Gamma \backslash G)^K$ given by
$$\lVert f \rVert_{2\ell}^2 \;\; = \;\; \lVert f \rVert^2 \; + \; \lVert (1-\Omega) \, f \rVert^2 \; + \; \lVert (1-\Omega)^2 \, f \rVert^2 \; + \; \dots \; + \; \lVert (1-\Omega)^{\ell} \, f \rVert^2$$
where $\lVert \, \cdot \, \rVert$ is the usual norm on $L^2(\Gamma \backslash G)$, induces a topology on $C^{\infty}_c(\Gamma \backslash G)^{K}$ that is equivalent to the topology induced by the family $\{\nu_{\alpha}: \alpha \in \mathcal{U}\mathfrak{g}^{\leq \, 2\ell}\}$ of seminorms and with respect to which $H^{2\ell}(X)$ is a Hilbert space.

\begin{proof} Let $\mathfrak{g} = \mathfrak{p} + \mathfrak{k}$ be the Cartan decomposition of $\mathfrak{g}$, and let $\{X_i\}$ be a basis for $\mathfrak{g}$.  Then $\Omega \, = \, \sum_{i } X_i \, X_i^{\ast}$.  Let $\Omega_{\mathfrak{p}}$ and $\Omega_{\mathfrak{k}}$ denote the subsums corresponding to the subspaces $\mathfrak{p}$ and $\mathfrak{k}$ respectively.  The set $\Sigma$ of possible $K$-types of $\gamma \, f$, for $\gamma \in \mathcal{U}\mathfrak{g}^{\leq \ell}$, is finite.  Let $\lambda_{\sigma}$ denote the $\Omega_{\mathfrak{k}}$-eigenvalue of a function of $K$-type $\sigma$, and let $C$ be a number greater than the maximum value of $\{\lambda_{\sigma} : \sigma \in \Sigma \}$.  
By the Poincar\'{e}-Birkhoff-Witt theorem we may assume $\alpha$ is a monomial of the form $\alpha \; = \; x_1 \dots x_n \; y_1 \dots y_m$ where $x_i \in \mathfrak{p}$ and $y_i \in \mathfrak{k}$.  For $f \in C_c^{\infty}(\Gamma \backslash G)^K$,
$$\nu_{\alpha} f \;\; = \;\; \langle \alpha f , \alpha f \rangle_{L^2(\Gamma \backslash G)} \;\; = \;\; \langle x_1 \dots x_n \; f, \; x_1 \dots x_n f \rangle_{L^2(\Gamma \backslash G)} \;\;\;\;\; (x_i \in \mathfrak{p})$$

\begin{lemma} For $\varphi \in C_c^{\infty}(\Gamma \backslash G)$ and $\alpha = x_1 \dots x_n$ a monomial in $\mathcal{U}\mathfrak{g}$ with $x_i \in \mathfrak{p}$,
$$\langle \alpha \, \varphi, \alpha \, \varphi \rangle \;\;  \leq \;\; \langle (-\Omega + C)^n \, \varphi, \, \varphi \rangle$$
where $\langle \, , \, \rangle$ is the usual inner product on $L^2(\Gamma \backslash G)$.

\begin{proof}
We proceed by induction on $n = \mathrm{deg} \, \alpha$.  For $n=1$, $\alpha = x \in \mathfrak{p}$.
Let $\{X_i\}$ be self-dual basis for $\mathfrak{p}$ such that $X_1 = x$.  Then, 
$$\langle x   \varphi, \, x  \varphi \rangle \;\; \leq \;\; \sum_i \langle X_i \,  \varphi, \, X_i \, \varphi \rangle \;\; = \;\; - \sum_i \langle X_i^2 \,  \varphi, \varphi \rangle $$
$$\;\;\;\;\;\;\;\;\;\;\;\; = \; \langle -\Omega_{\mathfrak{p}} \, \varphi,  \varphi \rangle \; = \;\;  \langle (-\Omega + \Omega_{\mathfrak{k}}) \,  \varphi,  \varphi \rangle  \;\; \leq \;\;  \langle (-\Omega + C) \, \varphi,  \varphi \rangle$$
For $n > 1$, write $\alpha = x \gamma$, where $x = x_1$ and  $\gamma = x_2 \dots x_n$.  Then, by the above argument,  $\langle x \, \gamma \varphi, \; x \, \gamma \varphi \rangle \;\; \leq \;\; \langle (-\Omega + C)\, \gamma \varphi, \, \gamma \varphi \rangle$.  Since the operator $-\Omega + C$ is a strictly positive symmetric unbounded operator, there is an everywhere defined inverse $R$, by Friedrichs \cite{friedrichs34, friedrichs35}, which is a positive symmetric \emph{bounded} operator, and which, by the spectral theory for bounded symmetric operators, has a positive symmetric square root $\sqrt{R}$ in the closure of the polynomial algebra $\C[R]$.  Thus  $-\Omega + C$  has a symmetric positive square root, namely $(1-\sqrt{R})$, commuting with all elements of $\mathcal{U}\mathfrak{g}$.  Thus, 
$$ \langle (-\Omega + C)\, \gamma \varphi, \, \gamma \varphi \rangle\;\; = \;\;  \langle \gamma \; \sqrt{-\Omega + C} \;  \varphi,  \;  \gamma \; \sqrt{-\Omega + C} \;  \varphi\rangle$$
By inductive hypothesis,
\begin{eqnarray*}
\langle \gamma \; \sqrt{-\Omega + C} \;  \varphi,  \;  \gamma \; \sqrt{-\Omega + C} \;  \varphi\rangle & \leq & \langle (-\Omega + C)^{n-1} \; \sqrt{-\Omega + C} \;  \varphi,  \;  \sqrt{-\Omega + C} \;  \varphi\rangle\\
& = &  \langle (-\Omega + C)^{n} \,  \varphi,  \;   \varphi\rangle
\end{eqnarray*}
This completes the proof of the lemma. \noqed \end{proof}
\end{lemma}

Thus, for any $\alpha \in \mathcal{U}\mathfrak{g}$, there is a constant $C$, possibly depending on the degree of $\alpha$, such that  $\nu_{\alpha}(f) \, \ll \, \langle (-\Omega + C)^{\mathrm{deg} \, \alpha} \, f , f \rangle$ for all $f \in C_c^{\infty}(\Gamma \backslash G)^K$.  In fact, for right $K$-invariant functions, $(-\Omega + C)^{\mathrm{deg} \, \alpha} \, f \; = \; (-\Omega_{\mathfrak{p}} + C)^{\mathrm{deg} \, \alpha} \, f $.  Since $\Omega_{\mathfrak{p}}$ is positive semi-definite, multiplying by a positive constant does not change the topology.  Thus, we may take $C=1$.  That is, the subfamily $\{\nu_{\alpha} : \alpha = (1-\Omega)^{k}, k \leq \ell \}$ of seminorms on $C_c^{\infty}(\Gamma \backslash G)^{K}$ dominates the family $\{\nu_{\alpha} : \alpha \in \mathcal{U}\mathfrak{g}^{\leq \, 2\ell} \}$ and thus induces an equivalent topology.  
\end{proof}
\end{prop}

It will be necessary to have another description of Sobolev spaces.  Let
$$W^{2, \ell}(\Gamma \backslash G) \;\; = \;\; \{f \in L^2(\Gamma \backslash G): \alpha \, f \in L^2(\Gamma \backslash G) \, \text{ for all } \alpha \in \mathcal{U}\mathfrak{g}^{\leq \ell} \}$$
where the action of $\mathcal{U}\mathfrak{g}$ on $L^2(\Gamma \backslash G)$ is by distributional differentiation.  Give $W^{2,\ell}(\Gamma \backslash G)$ the topology induced by the seminorms $\nu_{\alpha} f \; = \; \lVert \alpha \, f \rVert_{L^2(\Gamma \backslash G)}^2$,  $\alpha \in \mathcal{U}\mathfrak{g}^{\leq \ell}$, and let $W^{2, \ell}(X)$ be the subspace of right $K$-invariants. 

\begin{prop} \label{afc_large_small} These spaces are equal to the corresponding Sobolev spaces: 
$$W^{2, \ell}(\Gamma \backslash G) \; = \; H^{\ell}(\Gamma \backslash G) \;\;\;\;\; \text{ and } \;\;\;\;\; W^{2, \ell}(X)\; = \;H^{\ell}(X)$$

\begin{proof}
It suffices to show the density of test functions in $W^{2, \ell}(\Gamma \backslash G)$.  Since $G$ acts continuously on $W^{2,\ell}(\Gamma \backslash G)$ by right translation, \emph{mollifications} are dense in $W^{2,\ell}(\Gamma \backslash G)$; see \ref{gpints_moll}.  By Urysohn's Lemma, it suffices to consider mollifications of continuous, compactly supported functions.   Let $\eta \in C_c^{\infty}(G)$ and $f \in C_c^0(\Gamma \backslash G)$.  Then $\eta \cdot f$ is a smooth vector, and for all $\alpha \in \mathcal{U}\mathfrak{g}$, $\alpha \cdot (\eta \cdot f) \;\; = \;\; (L_{\alpha} \eta) \cdot f$.  For $X \in \mathfrak{g}$, the action on $\eta \cdot f$ as a vector is
$$\frac{\partial}{\partial t}\bigg|_{t=0} e^{tX} \cdot (\eta \cdot f) \;\; = \;\; \frac{\partial}{\partial t}\bigg|_{t=0} e^{tX}  \; \cdot \;\; \int_G \eta(g) \;\; g\cdot f \, dg \;\; = \;\; \frac{\partial}{\partial t}\bigg|_{t=0}  \;\; \int_G \eta(g) \;\; (g \, e^{tX} )\cdot f \, dg $$
Now using the fact that $f$ is a \emph{function} and the group action on $f$ is by \emph{translation},
$$\big( X \cdot (\eta \cdot f)\big)(h) \;\; = \;\; \frac{\partial}{\partial t}\bigg|_{t=0}  \;\; \int_G \eta(g) \;  f(he^{tX}g) \, dg \;\; = \;\;  \frac{\partial}{\partial t}\bigg|_{t=0}  \;\; (\eta \cdot f)(h e^{tX})$$
Thus the smoothness of $\eta \cdot f$ as a vector implies that it is a genuine smooth function.  The support of $\eta \cdot f$ is contained in the product of the compact supports of $\eta$ and $f$.  Since the product of two compact sets is again compact, $\eta \cdot f $ is compactly supported.
\end{proof}
\end{prop}

\begin{note}
By Proposition \ref{afc_casimir_norm}, $H^{2\ell}(X) = W^{2,2\ell}(X)$ is a Hilbert space with norm
$$\lVert f \rVert^2_{2\ell} \;\; = \;\; \lVert f \rVert^2 \; + \; \lVert (1-\Omega) \, f \rVert^2 \; + \; \dots \; + \; \lVert (1-\Omega)^{\ell} \, f \rVert^2$$
where $\lVert \, \cdot \, \rVert$ is the usual norm on $L^2(\Gamma \backslash G)$, and $(1-\Omega)^k \, f$ is a distributional derivative.
\end{note}

\subsection{Spectral transform, inversion, and differentiation on Sobolev spaces}

\begin{prop}\label{afc_diff_components} For $\ell \geq 0$, the Laplacian extends to a continuous linear map from $H^{2\ell +2 }(X)$ to $H^{2\ell}(X)$; the spectral transform extends to a map on $H^{2\ell}(X)$; and for every $f \in H^{2\ell+2}(X)$,  $\mathcal{F} \big( (1-\Delta) f \big) \, = \, (1-\lambda_{\xi}) \cdot \mathcal{F}$.
\begin{proof}
By the construction of the Sobolev topology, the Laplacian is a continuous linear map from  $C^{\infty}(\Gamma \backslash G) \cap H^{2\ell+2}(\Gamma \backslash G)$ to $C^{\infty}(\Gamma \backslash G) \cap H^{2\ell}(\Gamma \backslash G)$.  Since the Laplacian preserves right-$K$-invariance, it extends to a continuous linear map, also denoted $\Delta$, from $H^{2\ell +2}(X)$ to $H^{2\ell}(X)$.  The spectral transform, defined on $C^{\infty}_c(\Gamma \backslash G)^K$ by the integral transform, extends by continuity to $H^{2\ell}(X)$.  This extension agrees with the extension to $L^2(X)$ coming from Plancherel.  By integration by parts, $\mathcal{F}\big(\Delta \varphi\big) = \lambda_{\xi} \cdot \mathcal{F}\varphi$, for $\varphi \in C_c^{\infty}(\Gamma \backslash G)^K$, so by continuity, $\mathcal{F} \big((1- \Delta) f \big) \, = \, (1-\lambda_{\xi}) \cdot \mathcal{F}f$ for all $f \in H^{2\ell +2}(X)$.
\end{proof}
\end{prop}

Let $\mu$ be the multiplication map $\mu(v) (\xi) \; = \; (1-\lambda_{\xi}) \cdot v(\xi)$  on functions on $\Xi$.  For $\ell \in \Z$, the weighted $L^2$-spaces $V^{2\ell} \; = \; \{ v \text{ measurable } : \mu^{\ell}(v) \in L^2(\Xi) \}$ with norms
$$\lVert v \rVert_{V^{2\ell}} \;\; = \;\; \lVert \mu^{\ell}(v) \rVert_{L^2(\Xi)} \;\; = \;\; \int_{\Xi} (1-\lambda_{\xi})^{\ell} \; |v(\xi)|^2  \; d\xi$$
are Hilbert spaces with $V^{2\ell + 2} \subset V^{2\ell}$ for all $\ell$.  In fact, these are dense inclusions, since truncations are dense in all $V^{2\ell}$-spaces.  The multiplication map $\mu$ is a Hilbert space isomorphism $\mu: V^{2\ell + 2} \; \to \; V^{2\ell}$ since for $v \in V^{2\ell +2}$,
$$\lVert \mu(v) \rVert_{V^{2\ell}} \;\; = \;\; \lVert \mu^{\ell+1}(v) \rVert_{L^2(\Xi)} \;\; = \;\; \lVert v \rVert_{V^{2\ell+2}}$$
The negatively indexed spaces are the Hilbert space duals of their positively indexed counterparts, by integration.  The adjoints to inclusion maps are genuine inclusions, since $V^{2\ell+2} \hookrightarrow V^{2\ell}$ is dense for all $\ell \geq 0$, and the adjoint map $\mu^{\ast}:(V^{2\ell})^{\ast} \to (V^{2\ell +2})^{\ast}$ is the multiplication map $\mu:V^{-2\ell} \to V^{-2\ell -2}$.

\begin{prop}\label{afc_sp_trans_isom}  For $\ell \geq 0$, the spectral transform $\mathcal{F}$ is an isometric isomorphism $H^{2\ell}(X) \to V^{2\ell}$.

\begin{proof}
On compactly supported functions, the spectral transform $\mathcal{F}$ and its inverse $\mathcal{F}^{-1}$ are given by integrals, which are certainly continuous linear maps.  The Plancherel theorem extends $\mathcal{F}$ and $\mathcal{F}^{-1}$ to isometries between $L^2(X)$ and $L^2(\Xi)$.  Thus $\mathcal{F}$ on $H^{2\ell}(X)$ is a continuous linear $L^2$-isometry onto its image.    Let $f \in H^{2\ell}(X)$. By Proposition \ref{afc_large_small}, the distributional derivatives $(1-\Delta)^k \, f$ lie in $L^2(X)$ for all $k \leq \ell$.  By the Plancherel theorem and Proposition \ref{afc_diff_components},
$$\lVert (1-\Delta)^{\ell} \, f \rVert_{L^2(X)} \;\; = \;\; \lVert \mathcal{F}\big((1-\Delta)^{\ell} f\big)\rVert_{L^2(\Xi)} \;\;= \;\; \lVert(1- \lambda_{\xi})^{\ell} \cdot \mathcal{F}f \rVert_{L^2(\Xi)}$$
Thus $\mathcal{F}\big(H^{2\ell}(X)\big) \subset V^{2\ell}$.  The following claim shows that $\mathcal{F}^{-1}\big(V^{2\ell}\big) \subset H^{2\ell}(X)$, completing the proof of the proposition. \noqed \end{proof}

\begin{claim*} For $v \in V^{2\ell}$, the distributional derivatives $(1-\Delta)^k \, \mathcal{F}^{-1}v$ lie in $L^2(X)$, for all $0 \leq k \leq \ell$.
\begin{proof}
For test function $\varphi$, the Plancherel theorem implies
$$\big((1-\Delta) \, \mathcal{F}^{-1}v \big) \varphi \;\; = \;\; \mathcal{F}^{-1}v \, \big( (1-\Delta) \varphi \big) \;\; = \;\; v \big( \mathcal{F}\, (1-\Delta) \varphi \big) $$
By Proposition \ref{afc_diff_components} and the Plancherel theorem,
$$ v \big( \mathcal{F}\, (1-\Delta) \varphi \big) \;\; = \;\;  v \big((1-\lambda_{\xi})\cdot  \mathcal{F} \varphi \big) \;\; = \;\; \big((1-\lambda_{\xi})\cdot v\big) \, (\mathcal{F} \varphi) \;\; = \;\; \mathcal{F}^{-1} \big( (1-\lambda_{\xi}) \cdot v \big) \, \varphi$$
By induction, we have the following identity of distributions:  $(1-\Delta)^k \; \mathcal{F}^{-1} v \; = \;  \mathcal{F}^{-1} \big( (1-\lambda_{\xi})^k \, v \big)$.  Since $\mathcal{F}$ is an $L^2$-isometry and $(1-\lambda_{\xi})^{k} \, v \in L^2(\Xi)$ for all $0 \leq k \leq \ell$, $(1-\Delta)^k \; \mathcal{F}^{-1} v$ lies in $L^2(X)$ for $0 \leq k \leq \ell$.
\end{proof}
\end{claim*}
\end{prop}

\begin{note}\label{afc_spectral_defn}
This Hilbert space isomorphism $\mathcal{F}: H^{2\ell}(X) \to V^{2\ell}$ gives a \emph{spectral} characterization of the $2\ell^{\text{th}}$ Sobolev space, namely the preimage of $V^{2\ell}$ under $\mathcal{F}$.
$$H^{2\ell}(X) = \{ f \in L^2(X) : (1-\lambda_{\xi})^{\ell} \cdot \mathcal{F}f(\xi) \in L^2(\Xi) \}$$
\end{note}

\subsection{Negatively indexed Sobolev spaces and distributions}

The theory of negatively indexed Sobolev spaces allows the use of spectral theory for solving differential equations involving certain \emph{distributions.}

\begin{definition} 
For $\ell >0$, the Sobolev space $H^{-\ell}(X)$ is the Hilbert space dual of $H^{\ell}(X)$.
\end{definition}

\begin{note} Since the space of test functions is a dense subspace of $H^{\ell}(X)$ with $\ell >0$, dualizing gives an inclusion of $H^{-\ell}(X)$ into the space of distributions.  Similarly, the adjoints of the dense inclusions $H^{\ell} \hookrightarrow H^{\ell -1}$ are literal inclusions $H^{-\ell+1}(X) \hookrightarrow H^{-\ell}(X)$.  The self-duality of $H^0(X) = L^2(X)$ implies that $H^{\ell}(X) \hookrightarrow H^{\ell - 1}$ for all $\ell \in \Z$.
\end{note}

\begin{prop} \label{afc_ext_to_neg} The spectral transform extends to an isometric isomorphism on negatively indexed Sobolev spaces $\mathcal{F}:H^{-2\ell} \to V^{-2\ell}$, and $\mathcal{F} ((1-\Delta) \, u) \;\; = \;\; (1-\lambda_{\xi}) \cdot \mathcal{F} u$ for any $u \in H^{2\ell}(X)$, $\ell \in \Z$.

\begin{proof}  To simplify notation, for this proof let $H^{2\ell} = H^{2\ell}(X)$.   Propositions \ref{afc_diff_components} and \ref{afc_sp_trans_isom} give the result for positively indexed Sobolev spaces, expressed in the following commutative diagram,
\begin{equation*}
\xymatrix@C=3.5pc@R=3.5pc{
\dots \ar[r]^{(1-\Delta)} & H^4 \ar[r]^{(1-\Delta)} \ar[d]_{\mathcal{F}}^{\approx} & H^2 \ar[r]^{(1-\Delta)} \ar[d]_{\mathcal{F}}^{\approx} & H^0 \ar[d]_{\mathcal{F}}^{\approx}\\
\dots \ar[r]^{\mu} & V^4 \ar[r]^{\mu} & V^2 \ar[r]^{\mu} & V^0 }
\end{equation*}
where $\mu(v) (\xi) \;\; = \;\; (1-\lambda_{\xi}) \cdot v(\xi)$, as above.  Dualizing, we immediately have the commutativity of the adjoint diagram.
\begin{equation*}
\xymatrix@C=3.5pc@R=3.5pc{
(H^0)^{\ast} \ar[r]^{(1-\Delta)^{\ast}} & (H^2)^{\ast} \ar[r]^{(1-\Delta)^{\ast}}  & (H^4)^{\ast} \ar[r]^{(1-\Delta)^{\ast}} & \dots \\
(V^0)^{\ast} \ar[r]^{\mu^{\ast}} \ar[u]_{\approx}^{\mathcal{F}^{\ast}}& (V^{-2})^{\ast} \ar[r]^{\mu^{\ast}} \ar[u]_{\approx}^{\mathcal{F}^{\ast}}& (V^{-4})^{\ast} \ar[r]^{\mu^{\ast}} \ar[u]_{\approx}^{\mathcal{F}^{\ast}} & \dots }
\end{equation*}
The self-duality of $L^2$ and the Plancherel theorem allow the two diagrams to be connected.
\begin{equation*}
\xymatrix@C=3.5pc@R=3.5pc{
\dots \ar[r]^{(1-\Delta)} & H^4 \ar[r]^{(1-\Delta)} \ar[d]_{\mathcal{F}}^{\approx} & H^2 \ar[r]^{(1-\Delta)} \ar[d]_{\mathcal{F}}^{\approx} & H^0 \ar@/_1pc/[d]_{\mathcal{F}}^{\approx} \ar[r]^{(1-\Delta)^{\ast}} & H^{-2} \ar[r]^{(1-\Delta)^{\ast}}  & H^{-4} \ar[r]^{(1-\Delta)^{\ast}} & \dots\\
\dots \ar[r]^{\mu} & V^4 \ar[r]^{\mu} & V^2 \ar[r]^{\mu} & V^0  \ar[r]^{\mu} \ar@/_1pc/[u]^{\approx}_{\mathcal{F}^{-1}} & V^2 \ar[r]^{\mu} \ar[u]_{\approx}^{\mathcal{F}^{\ast}}& V^4 \ar[r]^{\mu} \ar[u]_{\approx}^{\mathcal{F}^{\ast}} & \dots }
\end{equation*}
Since $V^{2\ell + 2}$ is dense in $V^{2\ell}$ for all $\ell \in \Z$, and $H^{2\ell} \approx V^{2\ell}$ for all $\ell \in\Z$, $H^{2\ell + 2}$ is dense in $H^{2\ell}$ for all $\ell \in \Z$.  Thus test functions are dense in \emph{all} the Sobolev spaces.  The adjoint map $(1-\Delta)^{\ast}: H^{-2\ell} \to H^{-2\ell -2}$ is the continuous extension of $(1-\Delta)$ from the space of test functions, since, for a test function $\varphi$, identified with an element of $H^{-2\ell}$ by integration,
$$\big((1 - \Delta)^{\ast} \Lambda_{\varphi}\big)(f) \;\; = \;\; \Lambda_{\varphi}\big((1-\Delta)f\big) \;\; = \;\; \langle \varphi, (1-\Delta)f \rangle \;\; = \;\; \langle (1-\Delta) \varphi, f \rangle \;\; = \;\; \Lambda_{(1-\Delta) \varphi} (f)$$
for all $f$ in $H^{2\ell+2}$ by integration by parts, where $\Lambda_{\varphi}$ is the distribution associated with $\varphi$ by integration and $\langle \, , \, \rangle$ denotes the usual inner product on $L^2(X)$.  The map $\big(\mathcal{F}^{\ast}\big)^{-1}$ on $H^{-2\ell}$ is the continuous extension of $\mathcal{F}$ from the space of test functions, since for a test function $\varphi$, 
$$\big(\mathcal{F}^{\ast} \, \Lambda_{ \mathcal{F} \varphi}\big) (f ) \;\;= \;\; \Lambda_{\mathcal{F} \varphi}\big( \mathcal{F}f \big) \;\; = \;\; \langle \mathcal{F} \varphi , \mathcal{F} f \rangle_{V^{2\ell}} \;\; = \;\; \langle \varphi, f \rangle_{H^{2\ell}} \;\; = \;\; \Lambda_{\varphi}(f)$$
for all $f \in H^{2\ell}$.  Thus, the following diagram commutes.
\begin{equation*}
\xymatrix@C=3.5pc@R=3.5pc{
\dots \ar[r]^{(1-\Delta)} & H^4 \ar[r]^{(1-\Delta)} \ar[d]_{\mathcal{F}}^{\approx} & H^2 \ar[r]^{(1-\Delta)} \ar[d]_{\mathcal{F}}^{\approx} & H^0 \ar[d]_{\mathcal{F}}^{\approx} \ar[r]^{(1-\Delta)} & H^{-2} \ar[r]^{(1-\Delta)} \ar[d]_{\mathcal{F}}^{\approx} & H^{-4} \ar[r]^{(1-\Delta)} \ar[d]_{\mathcal{F}}^{\approx} & \dots\\
\dots \ar[r]^{\mu} & V^4 \ar[r]^{\mu} & V^2 \ar[r]^{\mu} & V^0  \ar[r]^{\mu} & V^{-2} \ar[r]^{\mu} & V^{-4} \ar[r]^{\mu}  & \dots }
\end{equation*}
In other words, the relation  $\mathcal{F} \big((1- \Delta) \, u \big) = (1-\lambda_{\xi}) \cdot \mathcal{F}u$ holds for any $u$ in a Sobolev space.
\end{proof}
\end{prop}

\begin{definition} \label{loc_sob_sp} For a smooth manifold $M$, the positively indexed local Sobolev spaces $H^{\ell}_{\text{loc}}(M)$ consist of functions $f$ on $M$ such that for all points $x \in M$, all open neighborhoods $U$ of $x$ small enough that there is a diffeomorphism $\Phi: U \to \R^n$ with $\Omega = \Phi(U)$ having compact closure, and all test functions $\varphi$ with support in $U$, the function  $(f \cdot \varphi) \circ \Phi^{-1} : \Omega \longrightarrow \C$ is in the Euclidean Sobolev space $H^{\ell}(\Omega)$.
\end{definition}

Recall the Sobolev embedding theorem for local Sobolev spaces.  For a smooth manifold $M$, $H^{\ell +k}_{\text{loc}}(M) \; \subset \; C^k(M) \;\;\;\;\; \text{for } \, \ell > \mathrm{dim}(M)/2$.

\begin{prop}\label{afc_global_sob_emb}  For $\ell \, > \, \mathrm{dim}(G/K)/2$,  $H^{\ell + k}(X) \; \subset \;  C^k(G/K)$.
\begin{proof}
Since positively indexed global Sobolev spaces $H^{\ell}(X)$ lie inside the corresponding local Sobolev spaces, $H^{\ell}(X) \subset H^{\ell}_{\text{loc}}(G/K)$, and by local Sobolev embedding, $H^{\ell}_{\text{loc}}(G/K)  \subset C^k(G/K)$.
\end{proof}
\end{prop}

This embedding of global Sobolev spaces into $C^k$-spaces is used to prove that the integral defining spectral inversion for test functions can be extended to sufficiently highly indexed Sobolev spaces, i.e. the abstract isometric isomorphism $\mathcal{F}^{-1} \circ \mathcal{F}: H^{\ell}(X) \to H^{\ell}(X)$ is given by an integral that is convergent uniformly pointwise, when $\ell > \mathrm{dim}(G/K)/2$, as follows.

\begin{prop}\label{afc_cvgc_sp_expn}  For $f \in H^{s}(X)$ , $s > k + \mathrm{dim}(G/K)/2$,
$$f \;\; = \;\; \int_{\Xi} \mathcal{F}f(\xi) \, \Phi_{\xi} \,  \, d\xi \;\;\;\;\; \text{in } H^s(X) \text{ and } C^k(X)$$
\begin{proof}
Let $\{\Xi_n\}$ be a nested family of compact sets in $\Xi$ whose union is all of $\Xi$, $\chi_n$ be the characteristic function of $\Xi_n$, and $f_n$ be given by the $C^{\infty}(X)$-valued Gelfand-Pettis integral (see \ref{gpints_moll})
$$f_{n} \;\; = \;\; \int_{\Xi} \chi_n(\xi) \, \mathcal{F}f(\xi) \, \Phi_{\xi}  \, d\xi$$
Since $ \chi_n(\xi) \, \mathcal{F}f(\xi) $ is compactly supported,  $f_n \, = \, \mathcal{F}^{-1} (\chi_n \cdot \mathcal{F}f)$.  Thus, by Proposition \ref{afc_sp_trans_isom},
$$\lVert f_n - f_m \rVert_{H^s(X)} \;\; = \;\; \lVert  (\chi_n -\chi_m) \cdot \mathcal{F}f \rVert_{V^{s}}$$
Since $\mathcal{F}f$ lies in $V^{s}$, these tails certainly approach zero as $n, m \to \infty$.  Similarly,
$$\lVert f_n - f \rVert_{H^s(X)} \;\; = \;\; \lVert  (\chi_n -1) \cdot \mathcal{F}f \rVert_{V^{s}} \;\; \longrightarrow \;\; 0 \;\;\;\;\;\;\;\; \text{ as } n \to \infty$$
By Proposition \ref{afc_global_sob_emb}, $f_n$ approaches $f$ in $C^k(X)$.
\end{proof}
\end{prop}

The embedding of global Sobolev spaces into $C^k$-spaces also implies that \emph{compactly supported} distributions lie in global Sobolev spaces, as follows.

\begin{prop} \label{afc_cs_distns_sob} Any compactly supported distribution on $X$ lies in a global archimedian spherical automorphic Sobolev space.  Specifically, a compactly supported distribution of order $k$ lies in $H^{-s}(X)$ for all $s \, > \, k + \mathrm{dim}(G/K)/2$.

\begin{proof}
A compactly supported distribution $u$ lies in $\big(C^{\infty}(G/K)\big)^{\ast}$.  Since compactly supported distributions are of finite order, $u$ extends continuously to $C^k(G/K)$ for some $k \geq 0$.  Using Proposition \ref{afc_global_sob_emb} and dualizing, $u$ lies in $H^{-(\ell+k)}(X)$, for $\ell \, > \, \mathrm{dim}(G/K)/2$.
\end{proof}
\end{prop}

\begin{note} In particular, this implies that the Dirac delta distribution at the base point $x_o = \Gamma \cdot 1 \cdot K$ in $\Gamma \backslash G/K$ lies in $H^{-\ell}(X)$ for all $\ell \, > \, \mathrm{dim}(G/K)/2$.\end{note}

\begin{prop}\label{afc_sph_trans_on_distns}  For a compactly supported distribution $u$ of order $k$,
$$\mathcal{F}u \;\; = \;\; u(\Phi_{\xi}) \;\;\;\;\;\;\;\; \text{in } V^{-s}\;\;\;\; \text{where } s > k + \mathrm{dim}(G/K)/2$$
\begin{proof}
A compactly supported distribution $u$ of order $k$ lies in $H^{-s}$ for any $s  >k + \mathrm{dim}(G/K)/2$.  Let $f$ be any element in $H^{s}(X)$.  Then,
$$\langle \mathcal{F}f, \mathcal{F}u \rangle_{V^{s} \times V^{-s}} \;\; = \;\; \langle f, u \rangle_{H^{s}\times V^{s}} \;\; = \;\; u(f)$$
Since the spectral expansion of $f$ converges to it in the $H^{s}(X)$ topology by Proposition \ref{afc_cvgc_sp_expn}, 
$$u(f) \;\; = \;\; u \bigg(\lim_n \; \int_{\Xi_n} \mathcal{F}f(\xi) \, \Phi_{\xi} \;  \, d\xi\bigg) \;\; = \;\; \lim_n \; u \bigg(\int_{\Xi_n} \mathcal{F}f(\xi) \, \Phi_{\xi}\;  \, d\xi\bigg)$$
Since the integral is a $C^{\infty}(X)$-valued Gelfand-Pettis integral (see \ref{gpints_moll}) and $u$ is an element of $(C^{\infty}(X))^{\ast}$,
$$u \bigg(\int_{\Xi_n} \mathcal{F}f(\xi) \, \Phi_{\xi} \;  \, d\xi\bigg) \;\; = \;\; \int_{\Xi_n} \mathcal{F}f(\xi) \, u(\Phi_{\xi}) \;  \, d\xi$$
The limit as $n\to \infty$ is finite, by comparison with the original expression which surely is finite, and thus
$$\langle \mathcal{F}f, \mathcal{F}u \rangle_{V^{s} \times V^{-s}} \;\; = \;\; \int_{\Xi} \mathcal{F}f(\xi) \, u(\Phi_{\xi}) \;  \, d\xi \;\; = \;\; \langle \mathcal{F}f, u(\Phi_{\xi}) \rangle_{V^s \times V^{-s}}$$
Thus, $\mathcal{F}u = u(\Phi_{\xi})$ as elements of $V^{-s}$
\end{proof}
\end{prop}

\begin{note} This implies that the spectral transform of the Dirac delta distribution is $\mathcal{F}\delta = \Phi_{\xi}(x_o)$.\end{note}

\subsection{Gelfand-Pettis integrals and mollification} \label{gpints_moll}

We describe the vector-valued (weak) integrals of Gelfand \cite{gelfand36} and Pettis \cite{pettis38} and summarize the key results; see \cite{garrett2011}.  For $X,\mu$ a measure space and $V$ a locally convex, quasi-complete topological vector space, a Gelfand-Pettis (or weak) integral is a vector-valued integral $C_c^o(X,V) \to V$ denoted $f \to I_f$ such that, for all $\alpha \in V^{\ast}$, 
$$\alpha(I_f) \; = \; \int_X \alpha \circ f \; d\mu$$
where this latter integral is the usual scalar-valued Lebesgue integral.  

\begin{note}
 Hilbert, Banach, Frechet, LF spaces, and their weak duals are locally convex, quasi-complete topological vector spaces; see \cite{garrett2011}.
\end{note}

\begin{theorem}\label{gp_ints}
\emph{(i)} Gelfand-Pettis integrals \emph{exist}, are \emph{unique}, and satisfy the following \emph{estimate}:
$$I_f \; \in \; \mu(\mathrm{spt}f) \cdot \big(\text{closure of compact hull of } f(X) \big)$$
\emph{(ii)} Any continuous linear operator between locally convex, quasi-complete topological vetor spaces $T : V \to W$ commutes with the Gelfand-Pettis integral: $T(I_f) \, = \, I_{Tf}$.
\end{theorem}

For a locally compact Hausdorff topological group $G$, with Haar measure $dg$, acting continuously on a locally convex, quasi-complete vector space $V$, the group algebra $C_c^o(G)$ acts on $V$ by \emph{averaging}:
$$\eta \cdot v \; = \; \int_G \eta(g) \, g \cdot v \, dg$$

\begin{theorem}\label{smooth_vectors_dense} \emph{(i)}
Let $G$ be a locally compact Hausdorff topological group acting continuously on a locally convex, quasi-complete vector space $V$.  Let $\{\psi_i\}$ be an approximate identity on $G$.  Then, for any $v \in V$, $\psi_i \cdot v \to v$ in the topology of $V$.

\emph{(ii)} If $G$ is a Lie group, the space $V^{\infty}$ of smooth vectors is dense in $V$, since for a smooth approximate identity $\{\eta_i\}$ on $G$, the mollifications $\eta_i \cdot v$ are smooth.  In particular, for $X \in \mathfrak{g}$,  $X \cdot (\eta \cdot v) \;\; = \;\; (L_X \eta) \cdot v$.
\end{theorem}

\begin{note} For a \emph{function space} $V$, the space of smooth vectors $V^{\infty}$  is \emph{not necessarily} the subspace of \emph{smooth functions} in $V$.  Thus Theorem \ref{smooth_vectors_dense} does \emph{not} prove the density of smooth \emph{functions} in $V$.\end{note}

\section{Lattice points and the automorphic spectrum: an exact formula}

Now we return to the problem posed at the beginning of the paper: expressing the relationship between the number of lattice points in an expanding region in a symmetric space to the automorphic spectrum in an exact formula.

\subsection{Spectral identity: expressions for the automorphic fundamental solution for $(\Delta - \lambda_z)^{\nu}$}

Let $G$ be a complex semi-simple Lie group with finite center and $K$ a maximal compact subgroup.  Let $G = NAK$, $\mathfrak{g} = \mathfrak{n} + \mathfrak{a} + \mathfrak{k}$ be corresponding Iwasawa decompositions.  Let $\Sigma$ denote the set of roots of $\mathfrak{g}$ with respect to $\mathfrak{a}$, let $\Sigma^+$ denote the subset of positive roots (for the ordering corresponding to $\mathfrak{n}$), and let $\rho = \tfrac{1}{2} \sum_{\alpha \in \Sigma^+} m_{\alpha} \alpha$, $m_{\alpha}$ denoting the multiplicity of $\alpha$.  Let $\mathfrak{a}_{\C}^{\ast}$ denote the set of complex-valued linear functions on $\mathfrak{a}$.  Let $\Gamma$ be a arithmetic subgroup.  Let $\Phi_{\xi}$ be a spectral basis for $L^2(\Gamma \backslash G)$ in the sense of \ref{sp_param}.  Consider the differential equation on the arithmetic quotient $X = \Gamma \backslash G/K$:
$$(\Delta - \lambda_z)^{\nu} \; v_z \;\; = \;\; \delta_{x_o}$$
where the Laplacian $\Delta$ is the image of the Casimir operator for $\mathfrak{g}$, $\lambda_z =  z^2 - |\rho|^2$ for a complex parameter $z$, $\nu$ is an integral power, and $\delta_{x_o} = \delta_{\Gamma \cdot 1\cdot K}$ is Dirac delta at the basepoint in $\Gamma \backslash G/K$.

\begin{prop}\label{afc_sp_expn_fund_soln}For integral $\nu \, > \, (\mathrm{dim} \, X)/2$, $v_z$ is a continuous right-$K$-invariant function on $\Gamma \backslash G$ with the following integral representation:
$$v_z(g) \;\; = \;\; \int_{\Xi} \;\frac{\overline{\Phi}_{\xi}(x_o)}{(\lambda_{\xi} - \lambda_z)^{\nu}} \;\;  \Phi_{\xi}(g) \, d\xi$$

\begin{proof}
Since $\delta_{x_o}$ is a compactly supported distribution of order zero, by Proposition \ref{afc_cs_distns_sob}, it lies in the global automorphic Sobolev spaces $H^{-\ell}(X)$ for all $\ell > (\mathrm{dim} \, X)/2$.  Thus there is an element $v_z$ of $H^{-\ell + 2 \nu}(X)$ satisfying this equation.  The solution $v_z$ is unique in Sobolev spaces, since any $w_z$ satisfying
$$(\Delta - \lambda_z)^{\nu} \; w_z \;\; = \;\; \delta_{x_o}$$
must necessarily have the same spectral transform.  For $\nu > (\mathrm{dim} \, X)/2$, by Proposition \ref{afc_global_sob_emb}, the solution is continuous, and by Proposition \ref{afc_cvgc_sp_expn}, 
$$v_z(g) \;\; = \;\; \int_{\Xi} \mathcal{F} \, v_z(\xi) \, \Phi_{\xi}(g) \, d\xi\;\; = \;\; \int_{\Xi} \;\frac{\overline{\Phi}_{\xi}(x_o)}{(\lambda_{\xi} - \lambda_z)^{\nu}} \;\;  \Phi_{\xi}(g) \, d\xi$$
\end{proof}
\end{prop}

Let $u_z$ denote the solution to the differential equation on the free space $X = G/K$:
$$(\Delta - \lambda_z)^{\nu} \; u_z \;\; = \;\; \delta_{1 \cdot K}$$
This free space solution, described in the following theorem, is computed explicitly in \cite{decelles2011}.

\begin{theorem}\label{fmla_for_fundl_soln} When $G$ is of odd rank and $\nu = (n+1)/2 +d$, where $d$ is the number of positive roots, counted without multiplicities, and $n$ is the rank,
$$u_z(a) \;\; = \;\; C_G \; \cdot \; \prod_{\alpha\in \Sigma^+}  \, \frac{\alpha(\log a)}{2 \sinh(\tfrac{\alpha(\log a)}{2})} \; \cdot \;\; \frac{e^{-z|\log a|}}{z}  $$
where $C_G$ is an explicit constant depending on the group $G$.  When $G$ is of even rank and $\nu  =  (n/2) +d + 1$,
$$u_z(a) \;\; = \;\;  C_G \; \cdot  \; \prod_{\alpha\in \Sigma^+}  \, \frac{\alpha(\log a)}{2 \sinh(\tfrac{\alpha(\log a)}{2})} \;\cdot \; \frac{ |\log a|}{ z} \; \cdot \; K_1(z \, |\log a|)$$
where $K_1$ is the usual Bessel function.
\end{theorem}

\begin{prop} \label{poincare_series} For $\mathrm{Re}(z) \gg 1$, the Poincar\'{e} series
$$\operatorname{P\acute{e}}_z(g)\;\; = \;\;\sum_{\gamma \in \Gamma} u_z(\gamma \cdot g)$$
converges absolutely and uniformly on compacts to a continuous function on $\Gamma \backslash G/K$.  Moreover, it is of moderate growth, and it is square-integrable modulo $\Gamma$.


\begin{proof}
By Proposition \ref{cvgce_poincare}, it suffices to show that the free-space fundamental solution $u_z$ is of sufficient rapid decay.  Let $\lVert \, \cdot \,  \rVert$ be the gauge on $G$ with $\sigma_o>0$ such that 
$$\int_{G} \, \frac{1}{\lVert g\rVert^{\sigma}} \, dg \; < \; \infty \;\;\;\;\;\;\;\; (\text{for } \sigma > \sigma_o)$$
Note that the product over positive simple roots is bounded, so it suffices to show that there is a $\sigma >\sigma_o$ such that $|\log a | \cdot e^{-\mathrm{Re}(z) |\log a|} \, \ll \, \lVert a \rVert^{-\sigma}$.  We claim that $| \log_A(a) |$ and $\log \big(  \lVert a \rVert\big)$ are comparable.  On diagonal matrices $(a_i)$, the gauge is  $\lVert (a_i) \rVert \,= \, \max_{1 \leq i \leq n} \{a_i, a_i^{-1}\}$.  Taking logarithms,
$$\log \big(\lVert (a_i) \rVert \big) \;\; = \;\; \max_{1 \leq i \leq n} \{|\log a_i|\} \;\;\;\;\;\;\;\; (\ell^{\infty}\text{-norm on } \mathfrak{a})$$
On the other hand,
$$|\log_A a| \;\; =\;\; \bigg(\sum_{i=1}^n (\log a_i)^2 \bigg)^{1/2}  \;\;\;\;\;\;\;\; (\ell^{2}\text{-norm on } \mathfrak{a})$$
The usual comparison:
$$\frac{1}{\sqrt{n}} \cdot  \max_{1 \leq i \leq n} \{|\log a_i|\} \; \leq \; \bigg(\sum_{i=1}^n (\log a_i)^2 \bigg)^{1/2} \; \leq \;  \max_{1 \leq i \leq n} \{|\log a_i|\}$$
allows us to conclude that $u_z$ is of sufficient rapid decay, as follows:
$$|u_z(a)| \;\; \ll \;\; |\log a | \cdot e^{-\mathrm{Re}(z) |\log a|} \;\; =\;\; \frac{|\log a|}{\big(e^{|\log a|}\big)^{\mathrm{Re}(z)}} \;\; \ll \;\; \frac{\log\big(\lVert a \rVert\big)}{\lVert a\rVert^{\mathrm{Re}(z)}} \;\; \ll \;\; \frac{1}{\lVert a \rVert^{\mathrm{Re}(z) -1}}$$
\end{proof}
\end{prop}


\begin{note}For the case $G = SL_2(\C)$,
$$\mathrm{P\acute{e}}_z(1) \;\; = \;\; \sum_{\gamma \in \Gamma} \; \frac{r_{\gamma} \, e^{-(2z-1) r_{\gamma}}}{(2z-1) \, \sinh r_{\gamma}} $$
where $r_{\gamma}$ is the Cartan radius of $\gamma$.  Thus $e^{-(2z-1)r_{\gamma}}$ is the $-(2z-1)^{\text{th}}$ power of the length of the arc from the basepoint $x_o = 1\cdot K$ to its image $\gamma \cdot x_o$.  On the quotient, this arc becomes a closed geodesic, and the sum over $\Gamma$ is closely related to the Selberg zeta function associated to $\Gamma$. \end{note}

\begin{theorem}[Spectral identity] \label{spec_id} For $\mathrm{Re}(z) \gg 1$,
$$\mathrm{P\acute{e}}_z(g) \;\; =\;\; \int_{\Xi} \;\; \frac{\overline{\Phi}_{\xi}(x_o) \cdot \Phi_{\xi}(g)}{(\lambda_{\xi} - \lambda_z)^{\nu}}  \;\; d\xi \;\;\;\;\;\;\;\;(\text{uniformly pointwise})$$
where $u_z$ is the free-space fundamental solution in Theorem \ref{fmla_for_fundl_soln}.  In particular, when $G$ is of odd rank,
$$ \sum_{\gamma \in \Gamma} \;\; C_G \; \cdot \prod_{\alpha\in \Sigma^+}  \, \frac{\alpha(H(\gamma \cdot g))}{2 \sinh(\tfrac{\alpha(H(\gamma \cdot g))}{2})} \; \cdot \;\; \frac{e^{-z|H(\gamma \cdot g)|}}{z} \;\;\; = \;\;\;  \int_{\Xi} \;\; \frac{\overline{\Phi}_{\xi}(x_o) \cdot \Phi_{\xi}(g)}{(\lambda_{\xi} - \lambda_z)^{\nu}}  \;\; d\xi $$
and when $G$ is of even rank,
$$\sum_{\gamma \in \Gamma} \; C_G  \cdot  \prod_{\alpha\in \Sigma^+}  \, \frac{\alpha(H(\gamma \cdot g))}{2 \sinh(\tfrac{\alpha(H(\gamma \cdot g))}{2})} \;\cdot \; \frac{ |H(\gamma \cdot g)|}{ z} \; \cdot \; K_1(z \, |H(\gamma \cdot g)|) \;\;\; = \;\;\;  \int_{\Xi} \;\; \frac{\overline{\Phi}_{\xi}(x_o) \cdot \Phi_{\xi}(g)}{(\lambda_{\xi} - \lambda_z)^{\nu}}  \;\; d\xi$$
where $H(g)$ is defined by $g = k \cdot exp(H(g)) \cdot k'$, for $k, k' \in K$.

\begin{proof}
By Proposition \ref{poincare_series}, the Poincar\'{e} series $\mathrm{P\acute{e}}_z$ is an automorphic fundamental solution for $(\Delta - \lambda_z)^{\nu}$ in $L^2(\Gamma \backslash G /K) = H^0(\Gamma \backslash G/K)$.  By the uniqueness of solutions in Sobolev spaces, $\mathrm{P\acute{e}}_z = v_z$ in an $L^2$-sense.  But both functions are continuous, so by the uniqueness of continuous functions in an $L^2$-equivalence class, $\mathrm{P\acute{e}}_z = v_z$ in $C^0(\Gamma \backslash G/K)$.
\end{proof}
\end{theorem}


We refer to the sum over $\Gamma$ as the \emph{geometric expression} of the automorphic fundamental solution and the integral over $\Xi$ as the \emph{spectral expression} of the same.

\subsection{Explicit formula for smoothed lattice-point counting}

To extract tangible information from this spectral identity we use the following well-known variant on the classical Perron method, which has the virtue of absolute convergence.

\begin{lemma} \label{var_perron}
$$ \frac{1}{2\pi i} \int_{\sigma - iT}^{\sigma + i T} \frac{e^{sX}}{s(s+\theta)(s + 2\theta) \dots (s+\ell \theta)} \; ds \;\; = \;\; \begin{cases} \frac{(1-e^{-\theta X})^{\ell}}{\ell! \, \theta^{\ell}} + O_{\sigma}\big(\frac{e^{\sigma X}}{T^{\ell+1} \cdot X} \big) & \text{ if } X >0 \\O_{\sigma}\big(\frac{e^{\sigma X}}{T^{\ell+1} \cdot |X|} \big) & \text{ if } X<0 \end{cases}$$
\end{lemma}

\begin{theorem} For a complex semi-simple Lie group $G$ of odd rank, with maximal compact $K$ and a co-compact lattice $\Gamma$, the number of lattice points within an expanding region of the basepoint $x_o = 1 \cdot K$ is related to the automorphic spectrum by the following explicit formula:
\begin{eqnarray*}
\lefteqn{ \tilde{C}_G \; \cdot \sum_{\gamma: \; |\log a_{\gamma}|<X} \;\;  \prod_{\alpha \in \Sigma^+} \frac{\alpha(\log a_{\gamma})}{2\sinh \big(\tfrac{\alpha(\log a_{\gamma})}{2}\big)} \; \cdot \;   \frac{1}{\ell! \, \theta^{\ell}} \big( 1 \; - \;  e^{- \theta (X-|\log a_{\gamma}|)} \big)^{\ell}}\\
& = & \big(A^{\ell,\theta}_{|\rho|} (X) \cdot e^{|\rho| X} \;\; + \;\; B^{\ell,\theta}_{|\rho|}(X) \cdot e^{- |\rho| X}\big) \cdot |\Phi_1(x_o)|^2\\
& & \;\;\;\;\;\;\; + \;\; \sum_{\Xi - \{\Phi_1\}} |\Phi_{\xi}(x_o)|^2 \cdot \big( A^{\ell,\theta}_{z_{\xi}}(X) \, e^{z_{\xi}X}  \; + \; B^{\ell,\theta}_{z_{\xi}} \, e^{-z_{\xi} X} \; + \; \mathrm{Per}^{\ell, \theta}_{z_{\xi}}(X)\big)
\end{eqnarray*}
where $\tilde{C}_G$ is an explicit constant depending only on the group, $\Phi_1$ is the constant automorphic form, $z_{\xi}$ is given by $\lambda_{\xi} =  z^2_{\xi} - |\rho|^2$, $A^{\ell, \theta}_{z_{\xi}}(X)$ and $B^{\ell, \theta}_{z_{\xi}}(X)$ are polynomial in $X$, of degree $(\nu -1)$, and rational in $z_{\xi}$, and $\mathrm{Per}^{\ell, \theta}_{z_{\xi}}(X)$ is of exponential decay in $X$ and rational in $z_{\xi}$.

\begin{proof}
The compactness of $\Gamma \backslash G$ implies that the spectrum $\Xi$ is discrete, and for $G$ of odd rank, the spectral identity of Theorem \ref{spec_id}, evaluated at the basepoint, becomes
$$C_G \;\cdot \;\;\sum_{\gamma \in \Gamma} \;\; \prod_{\alpha \in \Sigma^+} \frac{\alpha(\log a_{\gamma})}{2\sinh \big(\tfrac{\alpha(\log a_{\gamma})}{2}\big)} \; \cdot \;  \frac{e^{-z \, |\log a_{\gamma}|}}{z} \;\; \;\; = \;\;\;\; \sum_{\xi \in \Xi} \;\; \frac{|\Phi_{\xi}(x_o)|^2}{(\lambda_{\xi} - \lambda_z)^{\nu}}$$
where $C_G$ is an explicit constant  depending only on the group $G$, and $a_{\gamma}$ is given by $\gamma = k \cdot a_{\gamma} \cdot k'$.    We apply a Perron integral transform
$$P_{\ell, \theta}(f)(X) \;\; = \;\;  \frac{1}{2\pi i} \int_{\sigma + i\R} f(z) \; \cdot \;  \frac{z \cdot e^{zX}}{z(z+\theta)(z + 2\theta) \dots (z+\ell \theta)} \; dz$$
to both sides of the identity.   On the geometric side,
\begin{eqnarray*}
P_{\ell, \theta} \left(\frac{e^{z \, |\log a_{\gamma}|}}{z} \right) & = & \int_{\sigma + i\R} \frac{e^{z(X- |\log a_{\gamma}|)}}{z(z+\theta)(z + 2\theta) \dots (z+\ell \theta)} \; dz\\
& & \\
& = & \begin{cases} (1 - e^{-\theta(X - |\log a_{\gamma}|)})^{\ell}/(\ell! \, \theta^{\ell}) & \text{if } X \; > \; |\log a_{\gamma}| \\ 0 & \text{if } X \; < \; |\log a_{\gamma}| \end{cases}
\end{eqnarray*}
by Lemma \ref{var_perron}.  Thus
$$P_{\ell, \theta}\big(v_z(x_o)\big) \;\; = \;\;  C_G \; \cdot \sum_{\gamma: \; |\log a_{\gamma}|<X} \;\;  \prod_{\alpha \in \Sigma^+} \frac{\alpha(\log a_{\gamma})}{2\sinh \big(\tfrac{\alpha(\log a_{\gamma})}{2}\big)} \; \cdot \;   \frac{1}{\ell! \, \theta^{\ell}} \big( 1 \; - \;  e^{- \theta (X-|\log a_{\gamma}|)} \big)^{\ell}$$
On the spectral side, write $\lambda_{\xi} =  z_{\xi}^2 - |\rho|^2$.  Then $\lambda_{\xi} - \lambda_z = -(z- z_{\xi} )(z+ z_{\xi})$, and 
$$P_{\ell,\theta}\left( \frac{ |\Phi_{\xi}(x_o)|^2}{(\lambda_{\xi} - \lambda_z)^{\nu}} \right) \;\; = \;\;  \frac{1}{2\pi i} \int_{\sigma + i\R}  \frac{(-1)^{\nu} \,|\Phi_{\xi}(x_o)|^2}{(z-  z_{\xi})^{\nu}(z + z_{\xi})^{\nu}} \; \cdot \;  \frac{e^{zX}}{(z+\theta)(z + 2\theta) \dots (z+\ell \theta)} \; dz$$
Move $(-1)^{\nu} \,|\Phi_{\xi}(x_o)|^2$ outside the integral, and evaluate by residues.  The residues of the poles at the spectrum $z = \pm z_{\xi}$ are
$$A^{\ell,\theta}_{z_{\xi}}(X) \;\; = \;\; \frac{1}{(\nu-1)!} \; \lim_{z \to z_{\xi}} \; \frac{\partial^{\nu-1}}{\partial z^{\nu-1}}  \; \bigg((z+ z_{\xi})^{-\nu} \; \frac{e^{zX}}{(z+ \theta) \dots (z + \ell \theta)} \bigg)$$
$$B^{\ell,\theta}_{z_{\xi}}(X) \;\; = \;\; \frac{1}{(\nu-1)!} \; \lim_{z \to - z_{\xi}} \; \frac{\partial^{\nu-1}}{\partial z^{\nu-1}}  \; \bigg((z- z_{\xi})^{-\nu} \; \frac{e^{zX}}{(z+ \theta) \dots (z + \ell \theta)} \bigg)$$
Visibly, these are polynomial in $X$ and rational in $z_{\xi}$.  When $\Phi_{\xi} = \Phi_1$ is the constant automorphic form, $\lambda_{\xi} = 0 \Rightarrow z_{\xi} = \pm |\rho|$.  The sum of the residues of the simple poles at  $z = -m\theta$ is:
$$\mathrm{Per}^{\ell, \theta}_{z_{\xi}}(X) \;\; = \;\; \frac{1}{\theta^{\ell-1}} \; \sum_{m=1}^{\ell} \frac{(-1)^{m-1} \, e^{-m \theta X}}{(m-1)! \, (\ell-m)! \, (z_{\xi}^2 - m^2 \theta^2)^{\nu}}$$
Thus,
\begin{eqnarray*}
(-1)^{\nu} \, P_{\ell, \theta}\big(\mathrm{P\acute{e}}_z(x_o)\big) & = &  \big(A^{\ell,\theta}_{|\rho|} (X) \cdot e^{ |\rho| X} \;\; + \;\; B^{\ell,\theta}_{|\rho|}(X) \cdot e^{|\rho| X}\big) \cdot |\Phi_1(x_o)|^2 \\
& & \;\;\;\;\;\;\; + \;\; \sum_{\Xi - \{\Phi_1\}} |\Phi_{\xi}(x_o)|^2 \cdot \big( A^{\ell,\theta}_{z_{\xi}}(X) \, e^{z_{\xi}X}  \; + \; B^{\ell,\theta}_{z_{\xi}} \, e^{-z_{\xi} X} \; + \; \mathrm{Per}^{\ell, \theta}_{z_{\xi}}(X) \big)
\end{eqnarray*}
Since $v_z(x_o) = \mathrm{P\acute{e}}_z(x_o)$, by Theorem \ref{spec_id},  we have the desired equality, with $\tilde{C}_G = (-1)^{\nu} \cdot C_G$.
\end{proof}
\end{theorem}

%
%
%

\subsection{Gauges on groups and convergence of Poincar\'{e} series}

We recall some general facts about gauges on groups and convergence of Poincar\'{e} series.  See \cite{wallach1988} or Appendix 1 of \cite{diaconu-garrett09}.  For a countably-based, locally compact Hausdorff, unimodular group $G$ with compact subgroup $K$, a \emph{gauge} $g \to \lVert g \rVert$ is a continuous positive real-valued function on $G$ such that:

\begin{quote}
(1) $\lVert e \rVert = 1$, $\lVert g \rVert \geq 1$, and $\lVert g^{-1} \rVert = \lVert g \rVert$ \\
(2) \emph{Submultiplicativity}: $\lVert gh \rVert \leq \lVert g\rVert \cdot \lVert h \rVert$ \\
(3) $K$\emph{-invariance}: $\lVert k \cdot g \rVert = \lVert g \rVert = \lVert g \cdot k \rVert$\\
(4) \emph{Integrability}: for some $\sigma_o >0$,
$$\int_{G} \, \frac{1}{\lVert g\rVert^{\sigma}} \, dg \; < \; \infty \;\;\;\;\;\;\;\; (\text{for } \sigma > \sigma_o)$$
\end{quote}


General reductive groups have gauges, and on $GL_n$ they admit a particularly simple description, in terms of the operator norm:
$$\lVert g \rVert \;\; =\;\; \max\big(|\, g \, |_{\text{op}}, \, |g^{-1}|_{\text{op}} \big) \;\;\;\;\;\;\;\;\text{where} \;\;\;\;\;\;\;\;| \, g \, |_{\text{op}} \;\; =\;\; \sup_{|x| \leq 1} \, |g \cdot x|$$

\begin{definition} If convergent, the \emph{Poincar\'{e} series} associated to a function $f$ on $G$ is
$$\mathrm{P\acute{e}}_{f}(g) \;\; = \;\; \sum_{\gamma \in \Gamma} \, f(\gamma \cdot g) $$\end{definition}

\begin{prop} \label{cvgce_poincare} \emph{(i)} For a discrete subgroup $\Gamma \subset G$,
$$\sum_{\gamma \in \Gamma} \, \frac{1}{\lVert \gamma \rVert^{\sigma}} \; < \; \infty \;\;\;\;\;\;\;\; (\text{for } \sigma > \sigma_o)$$
where $\sigma_o$ is the power describing the integrability of the gauge.

\emph{(ii)} If there exists $\sigma > \sigma_o$, such that  $|f(g)| \, \ll \, \lVert g \rVert^{-\sigma} $, then the Poincar\'{e} series associated to $f$ converges absolutely and uniformly on compact sets.  Moreover,  $|\mathrm{P\acute{e}}_f(g) | \, \ll \, \lVert g \rVert^{\sigma}$ for all $\sigma > \sigma_o$.

\emph{(iii)} If there exists $\sigma > \sigma_o$ such that  $|f(g)| \;\; \ll \;\; \lVert g \rVert^{-2\sigma}$, then the Poincar\'{e} series associated to $f$ is square integrable modulo $\Gamma$, i.e.
$$\int_{\Gamma \backslash G} \; |\mathrm{P\acute{e}}_f(g)|^2 \, dg \;\; < \;\; \infty$$
\end{prop}

\bibliography{lattice_pts,sobolev_refs,harm_an_afms,gelfand_pettis}{}
\bibliographystyle{plain}

\end{document}